\documentclass[10pt]{article}  
\usepackage{theorem}                 
\usepackage{amssymb}

\newtheorem{theorem}{Theorem}
\newtheorem{corollary}{Corollary}
\newtheorem{proposition}{Proposition}
\newtheorem{lemma}{Lemma}

\theorembodyfont{\rmfamily}
\newtheorem{remark}{Remark}

\newenvironment{definition}
{\smallskip\noindent{\bf Definition\/}:}{\smallskip\par}

\newenvironment{proof}
{\noindent{\em Proof\/}.}{{ $\Box$}\smallskip\par}

\newcommand{\CC}{{\Bbb C}}

\newcommand{\PP}{{\Bbb P}}

\newcommand{\RR}{{\Bbb R}}
\newcommand{\ZZ}{{\Bbb Z}}

\newcommand{\eps}{\varepsilon}

\title{Mirror symmetry, Kobayashi's duality, and Saito's duality}
\author{W.Ebeling
\thanks{Partially supported by the DFG-programme ''Global methods in
complex geometry' (Eb 102/4--3).
 }
\\
\parbox{9cm}{\small
\begin{center} Universit\"{a}t Hannover, Institut f\"{u}r Mathematik \\
        Postfach 6009, D-30060 Hannover, Germany \\
        E-mail: ebeling\symbol{'100}math.uni-hannover.de
\end{center}}
}

\date{}

\begin{document}

\maketitle

\begin{abstract}
M.~Kobayashi introduced a notion of duality of weight systems. We tone this notion slightly down to a
notion called coupling. We show that coupling induces a relation
between the reduced zeta functions of the monodromy operators of the corresponding singularities
generalizing an observation of K.~Saito concerning Arnold's strange duality. We show that the weight
systems of the mirror symmetric pairs of M.~Reid's list of 95 families of Gorenstein K3 surfaces in
weighted projective 3-spaces are strongly coupled. This includes Arnold's strange duality where the
corresponding weight systems are strongly dual in Kobayashi's original sense. We show that the same
is true for the extension of Arnold's strange duality found by the author and C.~T.~C.~Wall. 
\end{abstract}

\section*{Introduction}
The mirror symmetry of Calabi-Yau threefolds has attracted the attention of
many physicists and mathematicians. One- and two-dimensional Calabi-Yau
varieties
are elliptic curves and K3 surfaces respectively. It is well-known that there
exist 3 families of weighted projective elliptic plane curves. The cones over these
curves are the simple-elliptic singularities of type $\widetilde{E}_6$,
$\widetilde{E}_7$, and $\widetilde{E}_8$ (see below). They are self-dual with
respect to mirror symmetry.

M.~Reid classified and
listed all families of K3 weighted projective hypersurfaces with Gorenstein
singularities (unpublished). It turned out that there are 95 such families. The
cones over these surfaces are called simple K3 hypersurface singularities.
These singularities were classified and thus Reid's list was rediscovered by
T.~Yonemura \cite{Yo}.
These surfaces include compactifications of the 14 exceptional unimodal
hypersurface
singularities of V.~I.~Arnold. It is well-known that the mirror symmetry
between the corresponding families of K3 weighted projective hypersurfaces
corresponds to Arnold's strange duality (see e.g.\ \cite{D2}). S.-M.~Belcastro
\cite{Be} determined for which of the 95 families the mirror symmetric
family is
again in Reid's list.

V.~V.~Batyrev \cite{Ba} showed that the mirror symmetry of Calabi-Yau
hypersurfaces in toric varieties is related to the polar duality between their
Newton polytopes. M.~Kobayashi \cite{Ko} discovered that Arnold's strange
duality corresponds to a duality of weight systems and this is related to
Batyrev's result. K.~Saito \cite{S1,S2} observed that Arnold's strange duality
corresponds to  a duality between the characteristic polynomials of the
monodromy operators of the corresponding dual singularities. In
\cite[Lect.~3, Problem~8.5]{Yu} it is
asked whether there are any possible relations among all these dualities
and mirror
symmetry. Here we give a partial answer to this question extending
\cite{E1} where it was
shown that Saito's duality can be derived from polar duality.

We consider weight systems $(a_1, \ldots, a_n;h)$ with
$$0 < a_0:=h - \sum_{i=0}^n a_i, \quad a_0 \vert h.$$
Let $f(x_1, \ldots , x_n)$ be polynomial of weighted degree $h$. Then the hypersurface 
$\widetilde{X}$ in the weighted projective space given by 
$$x_0^{h/a_0} + f(x_1, \ldots , x_n) =0$$
is a Calabi-Yau hypersurface (if it is quasismooth).
We introduce a notion of coupling of such weight systems which tones down Kobayashi's notion of
duality. We relate this to polar
duality in the same way as in \cite{E1}. The basic notion is the notion of
a weighted magic square $C$. 
The partner weight system corresponds to
the transpose of this matrix. 
The natural $\CC^\ast$-action on $\CC^n$ induces a monodromy transformation on the homology of the
fibre
$$F= \{ (x_1, \ldots , x_n) \in \CC^n \, | \, f(x_1, \ldots, x_n)=1 \}.$$
We consider the reduced zeta function $\widetilde{\zeta}_C(t)$
of this monodromy operator. We indicate how this rational
function can be computed from the matrix $C$. 
We show that the function $\widetilde{\zeta}_{C^t}(t)$ associated to the transpose matrix $C^t$ is in a
sense dual to $\widetilde{\zeta}_C(t)$ which generalizes Saito's duality and coincides
with it in the case when $n=3$ and $a_0=b_0=1$.

Then we investigate the coupling of weight systems for the weight systems of
Belcastro's list of mirror symmetric pairs inside Yonemura's list of 95
weight systems. It turns out that for any mirror symmetric pair the
corresponding weight systems are (strongly) coupled. The cases of Arnold's
strange duality are exactly those with $a_0=b_0=1$ where the Saito
duality holds in its strong form. Here the corresponding weight systems are strongly dual in
Kobayashi's original sense. The 31 cases with
$a_0=1$ and $f(x_1,x_2,x_3)=0$ having an isolated singularity at the origin are compactifications of 
the 31 Fuchsian singularities classified by
I.~Dolgachev \cite{D0}, I.~G.~Sherbak \cite{Sh}, and Ph.~Wagreich
\cite{Wag}. Many, but not all, of them have mirror symmetric partners inside
Yonemura's list. In \cite{E2} we asked whether the mirror symmetric
families to the Fuchsian singularities not involved in Arnold's strange
duality and its extension by the author and C.~T.~C.~Wall are realized by
singularities. Here we find for 7 of these Fuchsian singularities
singularities which are related to these singularities in a way explained in Section~3.

Finally we consider the extension of Arnold's strange duality found by the
author and C.~T.~C.~Wall. This again corresponds to mirror symmetry. Here
also weighted complete intersections in weighted projective 4-spaces are
involved. We associate weight systems to these varieties and we
show that the mirror symmetric pairs have (strongly) dual weight systems.

The author is grateful to N.~Yui for drawing his attention to
Belcastro's paper. He would like to thank Ch.~Okonek for useful discussions.

\section{Duality of weight systems}
An $(n+1)$-tuple of positive integers $W_{\bf a}=(a_1, \ldots , a_n;h)$
is called a {\em weight system}. The integers
$a_i$ are called the {\em weights} of $W_{\bf a}$ and $h$ is called the {\em
degree} of $W_{\bf a}$.

Two weight systems $W=(a_1, \ldots, a_n;h)$ and $W'=(a'_1, \ldots,
a'_n;h')$ are {\em
equivalent} if there exists a permutation $\sigma \in \frak{S}_{n}$ and a
rational
number $\lambda$ such that $\lambda a_{\sigma(i)} = a'_i$ for $i=1, \ldots
, n$ and
$\lambda h=h'$.
The weight system is called {\em reduced} if
$${\rm gcd}(a_1, \ldots , a_n)=1.$$

Each equivalence class contains a unique reduced weight system satisfying
$$a_1 \leq \ldots \leq a_n.$$ 
Let 
$$a_0:= h - \sum_{i=1}^n a_i.$$
In the sequel we shall assume that our weight system is reduced, satisfies  $a_1 \leq \ldots \leq
a_n$, and that $a_0 \neq 0$.

\begin{sloppypar}

If $a_0 >0$ and $a_0 \vert h$ then we shall call the weight system a {\em Calabi-Yau weight system}.
The reason for this is the following: 
Let $\PP(a_0,{\bf a})=\PP(a_0, \ldots , a_n)$ be the weighted complex projective space of
weight $(a_0, \ldots , a_n)$, i.e.\ the projective variety ${\rm Proj}\,
\CC[x_0, \ldots , x_n]$ where the degree of $x_i$ is $a_i$. Denote by $(x_0:
\ldots : x_n)$ the natural homogeneous coordinates of $\PP(a_0, {\bf a})$.
Let $f(x_1, \ldots, x_n)$ be an equation of weighted degree $h$ and define 
$$\widetilde{f}(x_0, x_1, \ldots, x_n) := x_0^{h/a_0} + f(x_1, \ldots , x_n).$$ 
Consider the hypersurface $\widetilde{X}:=\widetilde{f}^{-1}(0)$ in $\PP(a_0, {\bf a})$. Let $\CC^{n+1}$
be the affine
$(n+1)$-space with coordinates $(x_0,
\ldots , x_n)$. Assume that the hypersurface $\widetilde{X}$ is quasismooth, i.e.\ the
cone $C_{\widetilde{X}} = \{\widetilde{f}=0\}$ over $\widetilde{X}$ in $\CC^{n+1}$ is smooth outside
of the origin. By \cite[Theorem~3.3.4]{D1} the dualizing sheaf $\omega_{\widetilde{X}}$
satisfies $\omega_{\widetilde{X}}={\cal O}_{\widetilde{X}}$. Therefore $\widetilde{X}$ is a
(possibly singular) Calabi-Yau variety.

\end{sloppypar}

We recall some definitions of \cite{Ko}.
Let $W_{\bf a} = (a_1, \ldots , a_n;h)$ and $W_{\bf b} = (b_1, \ldots ,
b_n;k)$ be two
weight systems.

\begin{definition} Let $C$ be an $n \times n$ matrix with
entries in
the non-negative integers. The matrix $C$ is called a {\em weighted
magic square} of weight $(W_{\bf a},W_{\bf b})$ if
\begin{eqnarray*}
C (a_1, \ldots ,a_n)^t & = & (h, \ldots ,h)^t \mbox{ and}\\
(b_1, \ldots ,b_n) C & = & (k, \ldots, k).
\end{eqnarray*}
\end{definition}

Let $C=(c_{ij})$ be a weighted magic square
of weight $(W_{\bf a},W_{\bf b})$. 
Let $B$ be the
$n \times n$ matrix $(c_{ij}-1)$.  Let $A$ be the inverse matrix of $B$.
By \cite[Lemma~2.3.5(1)]{Ko}, $(\det C)/h=(\det B)/a_0$ and
$(\det C)/k=(\det B)/b_0$ and both numbers are integers.

\begin{lemma} \label{Lemma1}
We have
\begin{eqnarray*}
A (1, \ldots, 1)^t & = & ( {\textstyle \frac{a_1}{a_0}, \ldots , \frac{a_n}{a_0} })^t ,\\
(1, \ldots , 1) A & = & (  {\textstyle \frac{b_1}{b_0}, \ldots , \frac{b_n}{b_0} }).
\end{eqnarray*}
\end{lemma}

\begin{proof} By definition, $BA (1, \ldots , 1)^t=(1, \ldots , 1)^t$.
We have
\begin{eqnarray*}
B \left( \begin{array}{c} a_1 \\ \vdots \\a_n \end{array} \right)  & = &
\left( \begin{array}{ccc} c_{11}-1 & \ldots & c_{1n}-1 \\
 \vdots & \ddots & \vdots \\
c_{n1}-1 & \ldots & c_{nn}-1 \end{array} \right)
\left( \begin{array}{c} a_1 \\ \vdots \\a_n \end{array} \right) \\
& = &
\left( \begin{array}{c} h - \sum_{i=1}^n a_i \\ \vdots \\ h - \sum_{i=1}^n a_i
\end{array} \right)
= \left( \begin{array}{c} a_0 \\ \vdots \\a_0 \end{array} \right).
\end{eqnarray*}
This implies the first claim. The second claim follows in the same way.
\end{proof}

It follows from Lemma~\ref{Lemma1} that the weight systems $W_{\bf a}$ and $W_{\bf b}$ can be
retrieved from the matrix $C$.

\begin{definition}
A weighted magic square $C$ of weight $(W_{\bf a},W_{\bf b})$ is called {\em primitive} if $|\det
C|=h=k$.

The weight systems $W_{\bf a}$ and $W_{\bf b}$ are called {\em dual} if
there exists a primitive weighted magic square of weight $(W_{\bf a},W_{\bf b})$.

Two dual weight systems are called {\em strongly dual} if any row and any
column of
$C$ contains at least one zero.
\end{definition}

If two weight systems $W_{\bf a}$ and $W_{\bf b}$ are dual, then it follows that $k=h$ and
$b_0=a_0$.
We tone down this definition to include the case when $a_0 \neq b_0$. 

\begin{definition}
A weighted magic square $C$ of weight $(W_{\bf a},W_{\bf b})$ is called {\em almost primitive} if
$|\det C|=hb_0=ka_0$.

The weight systems $W_{\bf a}$ and $W_{\bf b}$ are called {\em coupled} if
there exists an almost primitive weighted magic square of weight $(W_{\bf a},W_{\bf b})$.

Two coupled weight systems are called {\em strongly coupled} if any row and any
column of
$C$ contains at least one zero.
\end{definition}

Let $C=(c_{ij})$ be a weighted magic square
of weight $(W_{\bf a},W_{\bf b})$. 
We now assume that $a_0 >0$. We show that the coupling of weight systems is related to
the  polar duality of associated Newton polytopes (cf.\ \cite{Ko,E1}). 

\begin{sloppypar}

\begin{definition} The $(n-1)$-simplex $\Gamma({\bf a})$ which is the
convex hull of
the row vectors of the matrix $C$ in $\RR^n$ is called a {\em Newton
diagram} of the weight system $W_{\bf a}$.

The $(n-1)$-simplex $\Delta({\bf a})$ which is the convex hull of the
vectors $(h/a_1, 0, \ldots , 0)$, \dots , $(0, \ldots , 0, h/a_n)$ in
$\RR^n$ is called the {\em full Newton diagram} of the weight system
$W_{\bf a}$.

Let $\overline{\Delta}({\bf a})$ be the $n$-simplex which is obtained from $\Delta({\bf a})$ 
by taking the
convex hull with the origin in $\RR^n$ and translating it by the vector $(-1,
\ldots ,-1)$, i.e.\ $\overline{\Delta}({\bf a})$ is the convex hull of the
vectors
$(-1+h/a_1,-1, \ldots ,-1)$, \dots , $(-1, \ldots ,-1, -1+h/a_n)$, $(-1,
\ldots ,-1)$.
\end{definition}

\end{sloppypar}

\begin{definition} Let $M \subset \RR^n$. Let $\langle \, ,\, \rangle$ denote
the Euclidean scalar product of $\RR^n$. The {\em polar dual} of $M$ is the
following subset of $\RR^n$:
$$M^\ast := \{ y \in \RR^n \, | \, \langle x,y \rangle \geq -1 \mbox{ for
all } x \in M \}.$$
\end{definition}

\begin{lemma} \label{Lemma2}
The polar dual $\overline{\Delta}({\bf a})^\ast$ is the $n$-simplex with
vertices
$v_1:=(1,0, \ldots ,0)$, \dots , $v_n:=(0, \ldots ,0,1)$, $v_0:=(-a_1/a_0,
\ldots , -a_n/a_0)$.
\end{lemma}

\begin{proof} \cite[Lemma~3.2]{Ko}
\end{proof}

\begin{proposition} \label{Prop1}
Let $\nabla$ be the convex hull of the vectors
$v_1-v_0$, \dots , $v_n-v_0$ in $\RR^n$. Then, in the coordinate system
given by taking the rows of $A$ as basis vectors, $\nabla$ is the convex
hull of the columns of $C$, hence a Newton diagram of the partner weight
system $W_{\bf b}$.
\end{proposition}

\begin{proof} By Lemma~\ref{Lemma2}, the claim is equivalent to the
following statement:
$$AC = \left( \begin{array}{cccc} 1+\frac{a_1}{a_0} & \frac{a_1}{a_0} &
\cdots & \frac{a_1}{a_0} \\
\frac{a_2}{a_0} & 1+ \frac{a_2}{a_0} & \cdots & \frac{a_2}{a_0} \\
\vdots & \vdots & \ddots & \vdots \\
\frac{a_n}{a_0} & \frac{a_n}{a_0} & \cdots & 1+ \frac{a_n}{a_0}
\end{array} \right) . $$
If $E$ denotes the $n \times n$
identity matrix and $\bf 1$ the matrix with all entries equal to 1,
then we have
$$AC = A(B + {\bf 1}) = AB + A {\bf 1} = E + A {\bf 1}.$$
Hence the claim follows from Lemma~\ref{Lemma1}.
\end{proof}

\section{Saito's duality}
Let $C=(c_{ij})$ be a weighted magic square
of weight $(W_{\bf a},W_{\bf b})$.
We shall associate a rational function $\widetilde\zeta_C(t)$ to the matrix $C$.

We consider the hypersurface $X$ in $\CC^n$
defined by the equation 
$f(x_1, \ldots , x_n)=0$, where
$$f(x_1, \ldots , x_n)=x_1^{c_{11}}x_2^{c_{12}} \cdots x_n^{c_{1n}} +
x_1^{c_{21}}x_2^{c_{22}} \cdots x_n^{c_{2n}} + \ldots + x_1^{c_{n1}}x_2^{c_{n2}} \ldots
x_n^{c_{nn}}.$$ 
Let 
$$F:= \{ (x_1, \ldots, x_n) \in \CC^n \, | \, f(x_1, \ldots , x_n)=1 \}$$
be the Milnor fibre of $f : (\CC^n,0) \to (\CC,0)$. 

If $W_{\bf a}$ is a Calabi-Yau weight system, then there is the following relation with the hypersurface
$\widetilde{X}$ in $\PP(a_0, {\bf a})$ defined by the equation
$$x_0^{h/a_0} +f(x_1, \ldots , x_n)=0.$$
Let $V$ be the hypersurface in $\PP({\bf a}):=\PP(a_1, \ldots , a_n)$ given by the equation $f(x_1, \ldots ,
x_n)=0$.
The mapping
$$\pi_{a_0} : \widetilde{X} \to \PP({\bf a}), \quad (x_0, x_1, \ldots, x_n) \mapsto (x_1, \ldots , x_n),$$
is a covering
of degree $h/a_0$ which is branched along the hypersurface $V$ and possibly along the singularities of $\PP({\bf
a})$. Let
$\widetilde{X}_0:= \widetilde{X} \setminus \pi_{a_0}^{-1}(V)$. Let $\widetilde{Y}$ be the hypersurface in
$\PP(1,{\bf a})$ given by the equation $x_0^h+f(x_1, \ldots , x_n)=0$. Then the mapping $\pi_1 : \widetilde{Y} \to
\PP({\bf a})$ is a covering of degree $h$ branched along the hypersurface $V$ and possibly along the
singularities of $\PP({\bf a})$. Then $\widetilde{Y}_0:= \widetilde{Y} \setminus \pi_1^{-1}(V)$ can be
identified with the Milnor fibre $F$ (cf.\ \cite{DD}). Therefore the
induced mapping $F=\widetilde{Y}_0 \to \widetilde{X}_0$ is a (possibly branched) covering of degree $a_0$. 

We have a $\CC^\ast$-action on $\CC^n$ defined by
$$\lambda \ast (x_1, \ldots , x_n) =(\lambda^{a_1}x_1, \ldots ,
\lambda^{a_n}x_n), \quad \lambda \in \CC^\ast.$$
Then the $\CC^\ast$-action induces a monodromy transformation $\theta : F \to F$
defined by
$$x \mapsto e^{2 \pi i/h} \ast x \quad (x \in F).$$
Let $\theta_\ast: \widetilde{H}_\ast(F) \to \widetilde{H}_\ast(F)$ be the induced operator
on the reduced homology of $F$. It is the
classical monodromy operator of the singularity $f(x_1, \ldots , x_n)=0$.
Let
$$\widetilde\zeta_C(t):=\prod_{p\ge0}
\left\{\det \left( \mbox{id} -t\cdot \theta_\ast\vert_{\widetilde
H_p(F)}\right)\right\}^{(-1)^p}$$ 
be the reduced zeta function of $\theta$. If $X$ has an isolated singularity at the
origin, the reduced zeta function is related to the characteristic polynomial $\phi_C(t)$
of the monodromy as follows:
$$\phi_C(t) = \left(\widetilde\zeta_C(t)\right)^{(-1)^{n-1}}.$$
The reduced zeta function $\widetilde\zeta_C(t)$ can be computed as follows (cf.\
\cite{EG2}). 

For $J \subset I_0= \{1, \ldots , n\}$ we denote by $\vert J
\vert$ the number of elements of $J$. For $J \neq \emptyset$, 
let $T_J := \{ x \in \CC^n \, | \, x_i = 0 \mbox{ for } i \not\in J, x_i \neq 0
\mbox{ for } i \in J \}$ be the (''coordinate'') complex torus of dimension
$\vert J \vert$, and
let $a_J := \gcd(a_j,\ j \in J)$. The integer $a_J$ is the order of the
isotropy group of the $\CC^\ast$-action on the torus $T_J$.
Let
$X_J:= X \cap T_J$, $F_J:= F \cap T_J$. The operator $\theta$ maps $F_J$ to itself; let
$\theta_J$ be the restriction of $\theta$ to $F_J$. We have
$$\widetilde\zeta_C(t)= (1-t)^{-1} \prod_{J:\vert J \vert \geq 1} \widetilde\zeta_{C,J}(t),$$
where $\widetilde\zeta_{C,J}(t)$ is the reduced zeta function of $\theta_J$.

Let $Z_J:= T_J/\CC^\ast$, $Y_J:= X_J/\CC^\ast$. Note that if $a_J$ does not divide $h$
then $Z_J
\setminus Y_J$ is empty. In this case, $\widetilde\zeta_{C,J}(t)=1$. Suppose $a_J
\vert h$. If we restrict the natural projection $T_J \setminus X_J \to Z_J
\setminus Y_J$ to $F_J$ then we get an $(h/a_J)$-fold covering $F_J \to Z_J \setminus Y_J$. The
transformation $\theta_J$ is a covering transformation of it and acts as  a cyclic
permutation of the $h/a_J$ points of a fibre. Therefore
$$\widetilde\zeta_{C,J}(t)=(1-t^{h/ a_J})^{ \chi(Z_J \setminus
Y_J)}$$ 
where $\chi(V)$ denotes the Euler characteristic of the topological space $V$.

The Euler characteristic $\chi(Z_J \setminus Y_J)$ can be computed as follows. A subset $J \subset
I_0$   is called {\em special} if there exists a subset $I \subset I_0$ with
$\vert I \vert = \vert J \vert$ such that $c_{ij}=0$ for
$i \in I$ and $j \not\in J$. Note that in particular $J=\emptyset$ and $J=I_0$ are special. For a
special subset
$J \neq \emptyset$ denote by $C_{IJ}$ the matrix $(c_{ij})^{i \in I}_{j \in J}$. Define
$C_{\emptyset \empty}:=(1)$. 

First assume that $\vert J \vert =1$.
Then
$Z_J= {\rm pt}$. The set
$Y_J$ is empty if and only if $J$ is special.  In this case, 
$$\widetilde\zeta_{C,J}(t)=(1-t^{h/a_J}).$$

Now suppose that $\vert J \vert \geq 2$. Then $\chi(Z_J) = 0$ and 
$\chi(Z_J \setminus Y_J) = - \chi(Y_J)$. Then
$Y_J \neq \emptyset$ if and only if $J$ is special or $J=I_0$. In this case, by \cite{BKKh, Kou} (see
also
\cite[(7.1)~Theorem]{Va}) we have
$$\chi(Y_J) = (-1)^{\vert J \vert} \frac{a_J \vert \det C_{IJ} \vert}{h}.$$
In particular, if $J=I_0$, then J is special, $a_J=1$ (since $W_{\bf a}$ is reduced), $C_{IJ}=C$, and
$$\chi(Y_J) = (-1)^n \frac{\vert \det C \vert}{h}.$$

For a  subset $J \subset I_0$ denote by $J'$ the complementary set $J':=I_0 \setminus J$. Note that
if $J$ is special for $C$ then $J'$ is special for $C^t$. For $J \neq \emptyset$ let $b_J :=
\gcd(b_j,\ j \in J)$. Define $a_\emptyset:=h$ and $b_\emptyset:=k$.

Summarizing we have proved the following theorem.

\begin{theorem} \label{ThmZeta}
The reduced zeta functions $\widetilde\zeta_C(t)$ and $\widetilde\zeta_{C^t}(t)$ can be computed
from the matrix $C$ as follows:
\begin{eqnarray*}
\widetilde\zeta_C(t) & = &  \prod_{J\, {\rm
special}} (1-t^{h/ a_J})^{(-1)^{\vert J \vert +1}a_J \vert \det C_{IJ} \vert /h},\\
\widetilde\zeta_{C^t}(t) & = &  \prod_{J\, {\rm
special}} (1-t^{k/ b_{J'}})^{(-1)^{\vert J' \vert +1} b_{J'} \vert \det C_{I'J'} \vert /k}.
\end{eqnarray*}
\end{theorem}

\begin{remark} \label{Rem1}
Let $X$ have an isolated singularity at the origin. Then its
Milnor number $\mu= \mbox{rank}\, H_{n-1}(F)$ is equal to
$$\mu=  (-1)^{n-1}\sum_{J\, {\rm special}} (-1)^{\vert J \vert +1}
\vert \det C_{IJ} \vert .$$ 
The dimension $\mu_0$ of the radical of $H_{n-1}(F)$ is equal to
$$\mu_0 = (-1)^{n-1}\sum_{J\, {\rm special}} (-1)^{\vert J
\vert +1} \frac{a_J\vert \det C_{IJ} \vert}{h}.$$ 
In addition, let $n=3$ and let $(a_1,a_2,a_3;h)$ be a Calabi-Yau weight system. By \cite[Theorem
3.3.4]{D1} the hypersurface $\widetilde{X}$ in
$\PP(a_0, {\bf a})$ is a simply-connected projective surface with dualizing sheaf
$\omega_{\widetilde{X}}= {\cal O}_{\widetilde{X}}$. Resolving its singularities (which are
rational double points) we get a non-singular K3 surface
with Picard number
$$\rho = 22 - (\mu - \mu_0).$$
If $\mu_0=0$, then by \cite[Proposition~1]{E0a} the discriminant of the Picard lattice,
i.e.\ the determinant of a matrix of the intersection form on the Picard group, is equal to
$$d=(-1)^{\rho -1} \widetilde\zeta_C(1) = (-1)^{\rho -1} \prod_{J\, {\rm
special}} \left(
\frac{h}{a_J}
\right)^{(-1)^{\vert J \vert +1} a_J\vert \det C_{IJ} \vert/h}.$$
\end{remark}

Following K.~Saito \cite{S1, S2}, for a rational function
$$
\psi(t) = \prod_{\ell\vert h}(1-t^\ell)^{\alpha_\ell}, \ \alpha_\ell\in\ZZ,
$$
we define the Saito dual (rational)
function $\psi^\ast(t)$ by 
$$
\psi^\ast(t)=\prod_{m\vert h}(1-t^m)^{-\alpha_{h/m}}.
$$
In particular, if $\sum_{\ell\vert h} \alpha_\ell =0$, then one has 
$$\psi^\ast(1) = \prod_{\ell\vert h} \left( \frac{h}{\ell} \right)^{-\alpha_\ell} = h^{\sum
\alpha_\ell} \prod_{\ell\vert h} \ell^{\alpha_\ell} = \psi(1).$$

\begin{corollary} \label{Cor1}
Let $C$ be primitive, $a_0=b_0=1$, and $n=3$. Then
$$\widetilde\zeta_{C^t}(t)=\widetilde\zeta_C^\ast(t).$$
\end{corollary}

\begin{proof} By the assumptions, we have $h=k$, $a_0=b_0=1$, and $I_0=\{ 1,2,3 \}$. We show that for
any special subset $J \subset I_0$ we have
$$a_J= \frac{h}{\vert \det C_{IJ} \vert}.$$
This is clear if $\vert J \vert =1$, $J=\emptyset$ or $J=I_0$. 

Therefore let $\vert J \vert =2$. For simplicity we
assume that $J=\{ 1,2 \}$. By Cramer's rule we have
$$a_1 = (c_{22}-c_{12}) \frac{h}{\det C_{IJ}}, \quad a_2 = (c_{11}-c_{21})
\frac{h}{\det C_{IJ}}.$$
This shows that $h/\vert \det C_{IJ} \vert$ divides $a_1$ and $a_2$ and hence $a_J$.
Let
$$a_J = e \frac{h}{\vert \det C_{IJ} \vert}$$
for some integer $e \geq 1$.
Then $e$ divides $c_{22}-c_{12}$ and $c_{11}-c_{21}$. If we subtract the second row of the matrix
$B$ from the first row then we obtain the matrix
$$\left( \begin{array}{ccc} c_{11}-c_{21} & c_{12}-c_{22} & 0 \\
                c_{21}-1 & c_{22} -1 & -1 \\
                c_{31}-1 & c_{32} -1 & c_{33} -1
\end{array} \right) .$$
Expanding the determinant of this matrix with respect to the first row we see that $e$ divides the
determinant of this matrix which is equal to $\det B =1$. This implies that $e=1$ and hence the
claim. 

Analogously, one can show that for $J$ special
$$b_{J'} = \frac{k}{\vert \det C_{I'J'} \vert}.$$

Hence it follows that for any special subset $J \subset I_0$
$$\frac{k}{b_{J'}} = \vert \det C_{I'J'} \vert = \frac{h}{\vert \det C_{IJ} \vert} = a_J.$$
Moreover,
$$(-1)^{\vert J' \vert +1} \frac{b_{J'} \vert \det C_{I'J'} \vert}{k} = (-1)^{\vert J' \vert +1} = -
(-1)^{\vert J \vert +1} \frac{a_J \vert \det C_{IJ} \vert}{h}.$$
Hence the claim follows from Theorem~\ref{ThmZeta}.
\end{proof}

\section{Simple K3 hypersurface singularities}
First consider the case $n=2$. Then the (Calabi-Yau) weight systems
corresponding to quasismooth plane curves are indicated in
Table~\ref{Tab1}. They are self-dual. The corresponding weighted magic
squares are given in that table. They are indicated as follows:
$$x^{c_{11}}y^{c_{12}}, x^{c_{21}}y^{c_{22}}.$$
The corresponding functions 
$f(x,y) = x^{c_{11}}y^{c_{12}} + x^{c_{21}}y^{c_{22}}$ have isolated singularities at the origin.
The characteristic polynomial
$\phi_C(t)$ of the monodromy operator satisfies $\phi_C^\ast(t)
=(\phi_C(t))^{-1}$ (cf. \cite{EG1}).
\begin{table}  \caption{Weighted elliptic plane
curves} \label{Tab1} $$
\begin{array}{|c|c|c|c|} \hline 
{\rm Name} & a_0,a_1,a_2;h & C & {\rm Dual} \\
\hline
\widetilde{E}_8 & 1,2,3;6 & x^3, y^2
& \widetilde{E}_8 \\ \widetilde{E}_7 & 1,1,2;4 & y^2, x^2y &
\widetilde{E}_7  \\ \widetilde{E}_6 & 1,1,1;3 & x^2y, xy^2 &
\widetilde{E}_6  \\ \hline \end{array} $$  \end{table}

Now consider the case $n=3$. Then the (Calabi-Yau) weight systems
corresponding to quasismooth surfaces have been classified by Reid
(unpublished) and Yonemura \cite{Yo}. The cones over these surfaces are called 
{\em simple K3 hypersurface singularities}.
Belcastro \cite{Be} determined the
mirror symmetric pairs inside that list.

\begin{theorem} \label{Thm1}
Let $W_{\bf a}$ and $W_{\bf b}$ be the weight systems of a mirror symmetric
pair of
simple K3 hypersurface singularities. Then $W_{\bf a}$ and $W_{\bf b}$ are
strongly coupled
weight systems.
\end{theorem}

For the proof of Theorem~\ref{Thm1} we indicate in each case an almost primitive
weighted magic square $C$ such that each row and each column of $C$ contains
at least one zero. This is done in Table~\ref{Tab2} for the weight systems
with
$a_i=b_j=1$ for some $i,j \in \{0,1\}$ and in Table~\ref{Tab3} for the
remaining cases.
We use the indexing of
\cite{Yo} for the weight systems. We list all the weight systems such that the
mirror family
is again in Yonemura's list. In the first column we indicate the index of
the weight
system. Let 
$$f(x,y,z)= x^{c_{11}}y^{c_{12}}z^{c_{13}} + x^{c_{21}}y^{c_{22}}z^{c_{23}}
+ x^{c_{31}}y^{c_{32}}z^{c_{33}}.$$
If $f(x,y,z)=0$ defines
an isolated hypersurface singularity
in Arnold's
\cite{Ar} or Wall's \cite{Wal} list of singularities, we give the name of
the singularity in the second column. In the case when $a_0=1$, $f(x,y,z)=0$ defines
a Fuchsian singularity (for the definition see \cite{E2}). In the cases where there is
a name missing we indicate the signature $\{g;\alpha_1, \ldots , \alpha_r\}$ of these singularities
(here 'nh' means that the central curve is non-hyperelliptic). In the third column we
list the weight system. In the 4th column we indicate the weighted magic
square $C$ in the following way:
$$x^{c_{11}}y^{c_{12}}z^{c_{13}}, x^{c_{21}}y^{c_{22}}z^{c_{23}}, x^{c_{31}}y^{c_{32}}z^{c_{33}}.$$
In the column preceding the
last one we indicate the index of the partner weight system. 

\begin{table}  \caption{Coupled weight systems with $a_i=b_i=1$ for $i=0$ or $i=1$} \label{Tab2}
$$
\begin{array}{|c|c|c|c|c|} \hline
{\rm No.} & {\rm Name} & a_0,a_1,a_2,a_3;h & C & {\rm Partner} \\
\hline
14 & E_{12} & 1,6,14,21;42 & x^7, y^3, z^2 & 14 \\
 & & 6,1,14,21;42 & x^{21}z,  y^3, z^2 & 28  \\
 & & 6,1,14,21;42 &  x^{28}y,  y^3, z^2 & 45   \\
 & & 14,1,6,21;42 &  x^{36}y, y^7,  z^2 & 51  \\
28 & & 3,1,7,10;21  &  x^{21}, y^3, xz^2 & 14  \\
 & & 3,1,7,10;21 &  x^{11}z,  y^3, xz^2 & 28  \\
 & & 3,1,7,10;21 &  x^{14}y, y^3, xz^2 & 45 \\
 & & 7,1,3,10;21 &  x^{18}y, y^7,  xz^2 & 51  \\
  45 & & 4,1,9,14;28 &  x^{28}, xy^3, z^2 & 14  \\
 & & 4,1,9,14;28 &  x^{14}z, xy^3, z^2 & 28  \\
 & & 4,1,9,14;28 &  x^{19}y, xy^3, z^2 & 45  \\
  & & 14,1,4,9;28 & x^{24}y, y^7, xz^3 & 51  \\
 51 & & 12,1,5,18;36 &  x^{36}, xy^7, z^2 & 14  \\
 & & 12,1,5,18;36 &  x^{18}z, xy^7, z^2 & 28  \\
  & & 18,1,5,12;36 &  x^{24}z, xy^7, z^3 & 45  \\
 & & 12,1,5,18;36  &  x^{31}y, xy^7, z^2 &  51  \\
 \hline
  50 & E_{13} & 1,4,10,15;30 &  x^5y, y^3, z^2 & 38   \\
& & 15,1,4,10;30 & x^{26}y, y^5z, z^3 & 77   \\
38 & Z_{11} & 1,6,8,15;30  &  x^5, xy^3, z^2 & 50 \\
 & & 15,1,6,8;30 & x^{22}z, y^5, yz^3 & 82   \\
77 & & 13,1,5,7;26  & x^{26}, xy^5, yz^3 & 50 \\
 & & 13,1,5,7;26 & x^{19}z, xy^5, yz^3 & 82 \\
82 & & 11,1,3,7;22  & x^{22}, y^5z, xz^3 & 38   \\
 & & 11,1,3,7;22 & x^{19}y, y^5z, xz^3 & 77 \\
 \hline
13 &  E_{14} & 1,3,8,12;24  & x^4z, y^3, z^2 & 20 \\
 & & 8,1,3,12;24 & x^{21}y, y^4z, z^2 & 59 \\
20 & Q_{10} & 1,6,8,9;24   &  x^4, y^3, xz^2 & 13 \\
 & & 8,1,6,9;24 & x^{15}z, y^4, yz^2 & 72   \\
59 & & 7,1,5,8;21 & x^{21}, xy^4, yz^2 & 13   \\
 & & 7,1,5,8;21 & x^{13}z, xy^4, yz^2 & 72 \\
 72 & & 5,1,2,7;15 & x^{15}, y^4z, xz^2 & 20  \\
 & & 5,1,2,7;15 & x^{13}y, y^4z, xz^2 & 59 \\
 \hline
78 & Z_{12} & 1,4,6,11;22  & x^4y, xy^3, z^2 & 78 \\ \hline
39 & Z_{13} & 1,3,5,9;18   &  x^3z, xy^3, z^2 & 60  \\
60 & Q_{11} & 1,4,6,7;18 & x^3y, y^3, xz^2 & 39 \\ \hline
22 & Q_{12} & 1,3,5,6;15   &  x^3z, y^3, xz^2 & 22 \\
 \hline
9  & W_{12} & 1,4,5,10;20  &  x^5, z^2, y^2z & 9 \\
 & & 4,1,5,10;20 & x^{15}y, z^2, y^2z & 71 \\
71 & & 3,1,4,7;15 & x^{15}, xz^2, y^2z & 9   \\
 & & 3,1,4,7;15 & x^{11}y, xz^2, y^2z & 71 \\
 \hline
37 & W_{13} & 1,3,4,8;16 &  x^4y, z^2, y^2z & 58 \\
 58 & S_{11} & 1,4,5,6;16 &  x^4, xz^2, y^2z & 37 \\ \hline
87 & S_{12} & 1,3,4,5;13   &  x^3y, xz^2, y^2z & 87 \\ \hline
4 & U_{12} & 1,3,4,4;12 &  x^4, y^3, z^3 & 4  \\ \hline
\end{array} $$
\end{table}

\begin{table} \caption{Coupled weight systems: remaining cases} \label{Tab3}
$$
\begin{array}{|c|c|c|c|c|}
\hline
{\rm No.} & {\rm Name} & a_0,a_1,a_2,a_3;h & C & {\rm Partner} \\
\hline
12 &  & 6,1,2,9;18 & x^9z, x^2y^8, z^2 & 27   \\
   &         & 6,1,2,9;18 & x^{16}y, x^2y^8, z^2 & 49 \\
27 &  & 8,2,3,11;24 & x^9y^2, y^8, xz^2 & 12   \\
49 &  & 14,2,5,21;42 & x^{16}y^2, xy^8, z^2 & 12  \\
\hline
40 &  & 7,1,2,4;14 & x^{10}z, x^2y^6, yz^3 & 81  \\
81 &  & 13,2,3,8;26 & x^{10}y^2, y^6z, xz^3 & 40 
\\ \hline
24 &  & 4,1,2,5;12 & x^{12}, x^2y^5, yz^2 & 11  \\ 
11 &  & 10,2,3,15;30 &  x^{12}y^2, y^5z, z^2 & 24   \\
\hline
6  &  & 2,1,2,5;10 &  x^5z, x^2y^4, z^2 & 26 \\
   &            & 2,1,2,5;10 &  x^8y, x^2y^4, z^2 & 34  \\
   &            & 5,1,2,2;10 &  x^8z, x^2y^4, yz^4 & 76  \\
26 &  & 4,2,5,9;20 &  x^5y^2, y^4, xz^2 & 6 \\ 
34 &  & 6,2,7,15;30 &  x^8y^2, xy^4, z^2 & 6  \\
76 &            & 13,2,5,6;26 &  x^8y^2, y^4z, xz^4 & 6 
\\ \hline
10 &  & 4,1,1,6;12 &  x^{11}y, y^6z, z^2 & 65   \\
& & 6,1,1,4;12  &  x^{11}y, y^8z, z^3 & 80  \\
 &         & 4,1,1,6;12 &  x^{12}, xy^{11}, z^2 & 46  \\
65 &  & 11,3,5,14;33 &  x^{11}, xy^6, yz^2 & 10  \\
80 &  & 22,4,5,13;44 &  x^{11}, xy^8, yz^3 & 10   \\
46 &         & 22,5,6,33;66 &  x^{12}y, y^{11}, z^2 & 10  \\
\hline
42 & Z_{2,0} & 1,1,3,5;10   &  x^5z, xy^3, z^2 & 68   \\
 & & 5,1,1,3;10 &  x^9y, y^7z, xz^3 & 92  \\
  & & 5,1,1,3;10 & x^9y, y^{10}, xz^3 & 83 \\
68 & Q_{17}  & 3,4,10,13;30 &  x^5y, y^3, xz^2 & 42 \\
92 & & 19,3,5,11;38 &  x^9z, xy^7, yz^3 & 42  \\
83 & & 27,4,5,18;54 &  x^9z, xy^{10}, z^3 & 42  \\
\hline
25 & Q_{3,0} & 1,1,3,4;9 &  x^6y, y^3, xz^2 & 43   \\
 & & 3,1,1,4;9    &  x^8y, y^5z, xz^2 & 88  \\
 & & 3,1,1,4;9 &  x^8y, y^9, xz^2 & 48   \\
43 & Z_{25}  & 4,3,11,18;36 &  x^6z, xy^3, z^2 & 25  \\
88 & & 9,2,5,11;27  &  x^8z, xy^5, yz^2 & 25  \\
48 &  & 16,3,5,24;48 &  x^8z, xy^9, z^2 & 25 \\
\hline
7 & X_{2,0} & 1,1,2,4;8    &  x^6y, y^2z, z^2 & 64  \\
64 & S_{17} & 3,4,7,10;24  &  x^6, y^2z, xz^2 & 7  \\
\hline
66 & S_{2,0}^\ast & 1,1,2,3;7 &  x^7, xz^2, y^2z & 35   \\
35 & W_{25} & 4,3,7,14;28 &  x^7y, z^2, y^2z & 66   \\ \hline
21 & 2;2 & 1,1,1,2;5 &  x^4y, y^3z, xz^2 & 86  \\
 &         & 1,1,1,2;5 &  x^4y, y^5, xz^2 & 30   \\
86 & V^\sharp NC^1_{18} & 5,4,7,9;25 &  x^4z, xy^3, yz^2 & 21 \\
30 & N_{33} & 8,5,7,20;40  &  x^4z, xy^5, z^2 & 21   \\
\hline
5  & 2;      & 1,1,1,3;6    &  x^5y, y^3z, z^2 & 56   \\
&         & 3,1,1,1;6 &  x^6, xy^5, yz^5 & 73  \\
56 & VNC^1_{13} & 5,6,8,11;30  &  x^5, xy^3, yz^2 & 5  \\
73 & & 25,7,8,10;50 &  x^6y, y^5z, z^5 & 5  \\ \hline
1  & 3;({\rm nh}) & 1,1,1,1;4 &  x^4, xy^3, yz^3 & 52  \\
52 & V'_{29} & 9,7,8,12;36 & x^4y, y^3z, z^3 & 1   \\ \hline
32 & & 2,2,3,7;14 & xy^4, x^4y^2, z^2 & 32 \\
\hline \end{array} $$
\end{table}

There are examples of strongly dual weight systems where the corresponding
families
of K3 surfaces are not mirror symmetric, e.g. (cf. \cite{Ko})
$$
\begin{array}{|c|c|c|c|c|} \hline
8 &  & 1,2,3,6;12 &  y^2z, x^3y^2, z^2 & 24  \\
24 &  & 1,2,4,5;12 & y^3, x^2y^2, xz^2 & 8  \\
\hline
\end{array}
$$
The weight systems correspond to the singularities $W_{1,0}$ and $Q_{2,0}$ respectively,
but the equations $f(x,y,z)=0$ are not equations of these singularities, they even
have non-isolated singularities at the origin. 

By inspection of the Tables 2 and 3, we see that the cases
$a_0=b_0=1$ are exactly the cases of
Arnold's strange duality. In these cases the
matrices $C$ are primitive and hence the corresponding weight systems are strongly
dual. In all other cases the weight systems are not strongly dual but only strongly coupled. 

The remaining singularities with $a_0=1$ are Fuchsian singularities of
signature $\{g;\alpha_1, \ldots , \alpha_r\}$ with $g >0$. They are coupled to weight systems which again correspond to
isolated singularities. We list these singularities together with their partners in
Table~\ref{TabFuchs}. Here $\rho$ denotes the Picard number of the K3 surface corresponding to the weight system on
the left-hand side as it can be found in the table of \cite{Be} and $d$ denotes the discriminant of the Picard
lattice. The numbers $\mu^\ast$, $\mu_0^\ast$, and $d^\ast$ are the Milnor number, the dimension of the radical, and
the discriminant of the Milnor lattice respectively of the singularity on the right-hand side. For the definition
of $\nu^\ast$ see below. The
singularities $Q_{17}$ and $S_{17}$ are bimodal singularities belonging to the list of Arnold \cite{Ar} of the 14 bimodal
exceptional singularities. The weight system of the singularity $V^\sharp
NC^1_{18}$ appears in the list of \cite[Appendix~1]{S1} of regular weight
systems with $\mu=24$ (see also \cite[Table~3]{E2}). The singularity
$VNC^1_{13}$ is a (minimally) elliptic hypersurface singularity and appears
in \cite[Table~2]{E0} (there we used the name $V^\sharp NC_{(1)}$). The singularities $Z_{25}$ and $W_{25}$ also
appear in the lists of Arnold \cite{Ar}. They have modality 4. The singularities $V'_{29}$ and $N_{33}$ do not
occur in the lists of \cite{Ar} and \cite{Wal}. Here we use as names the name of the series (according to
\cite{Ar}) to which they belong indexed by the Milnor number.
\begin{table} \caption{Fuchsian singularities with $g>0$ and their partners} \label{TabFuchs}
$$
\begin{array}{|c|c|c|c|c|c|c||c|c|c|c|c|c|}
\hline
{\rm No.} & {\rm Name} & \mu & \mu_0 & \rho & d & b_0 & d^\ast & \mu_0^\ast & \mu^\ast & \nu^\ast & {\rm Name} &
{\rm No.}
\\
\hline
42 & Z_{2,0} & 21 & 2 & 3 & 2 & 3 & -6 & 0 & 17 & 0 & Q_{17} & 68 \\
7  & X_{2,0} & 21 & 2 & 3 & 4 & 3 & -12 & 0 & 17 & 0 & S_{17} & 64 \\
21 & 2;2 & 24 & 4 & 2 & -5 & 5 & 25 & 0 & 24 & 0 & V^\sharp NC^1_{18} & 86  \\
 &  & &  &  &  & 8 & -10 & 0 & 33 & 6 & N_{33} & 30 \\
5 & 2; & 25 & 4 & 1 & 2 & 5 & -10 & 0 & 19 & 0 & VNC^1_{13} & 56 \\ 
25 & Q_{3,0} & 20 & 2 & 4 & -3 & 4 & -6 & 0 & 25 & 2 & Z_{25} & 43\\ 
66 & S_{2,0}^\ast & 20 & 2 & 4 & -7 & 4 & -14 & 0 & 25 & 2 & W_{25} & 35\\
1 & 3;({\rm nh}) & 27 & 6 & 1 & 4 & 9 & -12 & 0 & 29 & 6 & V'_{29} & 52\\
\hline
\end{array}
$$
\end{table}

\begin{sloppypar}

There is the following relation between these singularities. The K3 surfaces corresponding to the weight systems on
the left-hand side are compactifications of the corresponding Fuchsian singularities.  Let $g(x,y,z)=0$ be the
equation of a singularity on the right-hand side. Let 
$W_{\bf b}=(b_1,b_2,b_3;k)$ be the weight system of this singularity. Let $\widetilde{X}^\ast$ be the hypersurface
in
$\PP(b_0,b_1,b_2,b_3)$ given by the equation $w^{k/b_0}+g(x,y,z)=0$. As in
Section~2 we consider the natural mapping $\pi_{b_0} : \widetilde{X}^\ast \to \PP(b_1,b_2,b_3)$. Let
$\widetilde{X}^\ast_0 :=
\widetilde{X}^\ast \setminus \pi^{-1}_{b_0}(V^\ast)$ where $V^\ast$ is the hypersurface in $\PP(b_1,b_2,b_3)$
defined by
$g(x,y,z)=0$ and let $F^\ast$ be the Milnor fibre of 
$g: (\CC^3,0) \to (\CC,0)$. Then we have a mapping  $F^\ast \to \widetilde{X}^\ast_0$ which is a (possibly branched)
covering of degree $b_0$. Denote by $\nu^\ast$ the total branching order of this covering. One can easily see that
this covering is either unbranched or branched along the singularity $(0:0:1) \not\in V^\ast$ of
$\PP(b_1,b_2,b_3)$ of branching order $\nu^\ast$. Then we obtain from the Riemann-Hurwitz formula in all cases
$$\mu^\ast+\nu^\ast+1=b_0(\rho+3).$$

\end{sloppypar}

\section{Extension of Arnold's strange duality}
The author and C.~T.~C.~Wall \cite{EW} have found an extension of Arnold's
strange duality embracing also isolated complete intersection
singularities. Such a singularity is defined by the germ of an analytic
mapping $(g,f): (\CC^4,0) \to (\CC^2,0)$. It is weighted homogeneous of
weights $q_1, q_2,q_3,q_4$ and degrees $d_1, d_2$ where we assume $d_1 \leq
d_2$ and where we have $1+q_1+q_2+q_3+q_4=d_1+d_2$. We consider the
compactification of such a singularity in the weighted projective space
$\PP(1,q_1,q_2,q_3,q_4)$ with coordinates $w,x,y,z,t$ given by the
equations
\begin{eqnarray*} g(x,y,z,t) & = & 0, \\ f(x,y,z,t) + w^{d_2} & = & 0.
\end{eqnarray*}
More precisely, this correspondence embraces the following singularities.
We use the notation of \cite{E1}.
\begin{itemize}
\item[(a)] Arnold's 14 exceptional unimodal hypersurface singularities.
\item[(b)] The six bimodal hypersurface singularities $J_{3,0}$ (12),
$Z_{1,0}$ (40), $Q_{2,0}$ (24), $W_{1,0}$ (8), $S_{1,0}$ (63), $U_{1,0}$
(18). The compactifications of these singularities occur in Yonemura's
list. The index is indicated in brackets. The first three of these
singularities already occurred in Table~\ref{Tab3}.
\end{itemize}
The remaining singularities are ICIS defined by the germ of an analytic
mapping $(g,f): (\CC^4,0) \to (\CC^2,0)$ as above. Here we distinguish
between three types:
\begin{itemize}
\item[(c)] The singularities $J'_9$, $J'_{10}$,
$J_{11}$, $K'_{10}$, $K'_{11}$, $J'_{2,0}$, $K'_{1,0}$ where
$g(x,y,z,t)=xt-y^2$, $f(x,y,z,t)=f'(x,y,t)+z^2$ for some $f': (\CC^3,0) \to
(\CC,0)$.
\item[(d)] The singularities $L_{10}$, $L_{11}$, $M_{11}$, $L_{1,0}$,
$M_{1,0}$ where $g(x,y,z,t)=xt-yz$.
\item[(e)] The ICIS
$I_{1,0}$ given by \begin{eqnarray*} g(x,y,z,t) & = & x^3- yt,\\
f(x,y,z,t) & = & (a+1)x^3+yz+z^2+zt, \quad a \neq 0,1.
\end{eqnarray*}
\end{itemize}

The correspondence between these singularities is indicated in
Table~\ref{Tab4}. The compactifications of all these singularities are K3
surfaces and the dual families are mirror symmetric.

\begin{table} \caption{Extension of Arnold's strange duality} \label{Tab4}
$$
\begin{array}{|c|c|c|c|c|}
\hline
{\rm No.} & {\rm Name} & a_0,a_1,a_2,a_3;h & C & {\rm Dual}  \\
\hline
12 & J_{3,0} & 1,2,6,9;18 & x^3y^2, y^3, z^2 & J'_9  \\
40 & Z_{1,0} & 1,2,4,7;14 &  x^3y^2, xy^3, z^2 & J'_{10}  \\
24 & Q_{2,0} & 1,2,4,5;12 &  x^2y^2, y^3, xz^2 & J'_{11}  \\
8  & W_{1,0} & 1,2,3,6;12 &  x^3z, z^2, y^2z & K'_{10}, L_{10}  \\
63 & S_{1,0} & 1,2,3,4;10 &  x^3z, xz^2, y^2z & K'_{11}, L_{11}  \\
18 & U_{1,0} & 1,2,3,3;9  &  x^3y, y^2z, yz^2 & M_{11}  \\
\hline
 & J'_9 & 1,6,2,9;18 &  x^3, x^2y^3, z^2 & 12  \\
 & J'_{10} & 1,4,2,7;14 &  x^3y, x^2y^3, z^2 & 40  \\
 & J'_{11} & 1,3,2,6;12 &  x^2z, x^2y^3, z^2 & 24  \\
 & K'_{10}, L_{10} & 1,4,1,6;12 &  x^3, z^2, xy^2z & 8  \\
 & K'_{11}, L_{11} & 1,3,1,5;10 & x^3y, z^2, xy^2z & 63  \\
 & M_{11} & 1,3,1,4;9 &  x^3, xy^2z, yz^2 & 18  \\
 \hline
  & J'_{2,0} & 1,2,2,5;10 &  x^3y^2, x^2y^3, z^2 & J'_{2,0}  \\
  & L_{1,0}, K'_{1,0} & 1,2,1,4;8 &  x^2z, z^2, xy^2z  &
L_{1,0}, K'_{1,0}  \\
  & M_{1,0} & 1,2,1,3;7 &  x^3y, xy^2z, yz^2 & M_{1,0}  \\
  & I_{1,0} & 1,2,3,0;6 &  x^3, y^2z^2, y^2z & I_{1,0}  \\
\hline
\end{array} $$ \end{table}

We also relate this correspondence to a duality of weight systems. For this
purpose, we associate a Calabi-Yau weight system to an ICIS as follows. In
the cases (c) and (d) we associate the weight system
$(1,q_1,q_2-q_1,q_3;d_2)$ to the singularity $(X,0)$. Since $d_1=q_1+q_4$,
this is a Calabi-Yau weight system.

In case (e) we associate the weight system $(1,q_1, q_2,q_3-q_2;d_2)$ to
the singularity $I_{1,0}$. Since $d_1=q_2+q_4$, we have
$1+q_1+q_2+q_3-q_2=d_2$. However, $a_3:=q_3-q_2=0$.

Then we have the following extension of Theorem~\ref{Thm1}:

\begin{theorem} \label{Thm2}
Let $W_{\bf a}$ and $W_{\bf b}$ be the weight systems of a mirror symmetric
pair of the above singularities. Then $W_{\bf a}$ and $W_{\bf b}$ are
strongly dual.
\end{theorem}

In each case, a primitive weighted magic square
$C$ is indicated in Table~\ref{Tab4}. For a singularity of
type (b), the matrix $C$ corresponds to some points of the
Newton diagram. The corresponding function $f(x,y,z)$ has a non-isolated singularity at the
origin and does not define the given
one. These points differ from the points given in
\cite{E1}. The points there correspond to non-primitive (and even not almost primitive) matrices $C$.

For the ICIS of types (c)
or (d), one obtains some points of a Newton diagram of the Laurent
polynomial associated to the singularity in \cite{E1} by subtracting the
second column of $C$ from the first one. In case (e), one has
to subtract the third column from the second one.

In all cases, the reduced zeta function $\widetilde\zeta_C(t)$ differs from the characteristic
polynomial $\Delta(t)$ of the monodromy of the corresponding singularity
only in the exponents: If
$$\Delta(t) = \prod_{m|h} (1-t^m)^{\chi_m}$$
then
$$\widetilde\zeta_C(t) = \prod_{m|h} (1-t^m)^{\eps_m}$$
where $\eps_m=-1,0,1$ if
$\chi_m<0$, $\chi_m=0$, $\chi_m>0$ respectively. From Corollary~\ref{Cor1}
we get Saito's duality of the characteristic polynomials of the monodromy
up to the absolute value of the exponents.

For a generalization of the construction of the polar dual in \S 1 for this
extension of Arnold's strange duality which precisely yields
Saito's duality  we refer to \cite{E1}.

By inspection of Table~\ref{Tab4} we observe the strange fact that the
weight systems associated to the ICIS with the exception of $I_{1,0}$ again occur
in Yonemura's list. However, comparing
\cite[Table~3]{Be} and \cite[Table~6]{E0a}, we see that the Picard lattices
of the corresponding K3 surfaces are different.
If we omit the zero in the weight system of $I_{1,0}$,
we obtain the weight system of $\widetilde{E}_8$.


\end{document}